\documentclass[11pt]{amsart} 
\usepackage{amssymb}
\usepackage{amsmath}
\usepackage{setspace}
\usepackage{amsthm}
\usepackage{epsfig,color,colordvi,wasysym}
\usepackage{tipa}
\usepackage{tikz}
\usepackage{verbatim}

\newtheorem{theorem}{Theorem}

\def\x{\mathbf{x}}

\begin{document}

\title{The Chromatic Polynomials of Signed Petersen Graphs} % and Small Complete Graphs

\author{Matthias Beck}
\address{Department of Mathematics, San Francisco State University, San Francisco, CA 94132, USA}
\email{mattbeck@sfsu.edu}

\author{Erika Meza}
\address{Department of Mathematics, Loyola Marymount University, Los Angeles, CA 90045, USA}
\email{erikameza9@aol.com}

\author{Bryan Nevarez}
\address{Department of Mathematics, Queens College, CUNY, Flushing, NY 11367, USA}
\email{smooth\underline{ }bry@hotmail.com}

\author{Alana Shine}
\address{Department of Computer Science, University of Southern California, Los Angeles, CA 90089, USA}
\email{ashine@usc.edu}

\author{Michael Young}
\address{Department of Mathematics, Iowa State University, Ames, IA 50011, USA}
\email{myoung@iastate.edu}

\begin{abstract}
Zaslavsky proved in 2012 that, up to switching isomorphism, there are six different signed Petersen graphs and that they could be told apart by their chromatic polynomials, by showing that the latter give distinct results when evaluated at 3. He conjectured that the six different signed Petersen graphs also have distinct zero-free chromatic polynomials, and that both types of chromatic polynomials have distinct evaluations at \emph{any} positive integer.
We developed and executed a computer program (running in {\tt SAGE}) that efficiently determines the
number of proper $k$-colorings for a given signed graph; our computations for the signed Petersen
graphs confirm Zaslavsky's conjecture. We also computed the chromatic polynomials of all signed
complete graphs with up to five vertices.
\end{abstract}


% 05C22 Signed and weighted graphs
% 05A15 Exact enumeration problems, generating functions
% 05C15 Coloring of graphs and hypergraphs

\date{18 December 2014}

\thanks{We are grateful to Thomas Zaslavsky and an anonymous referee for comments on an earlier version of this paper, and we thank Ricardo Cortez and the staff at MSRI for creating an ideal research environment at MSRI-UP.
This research was partially supported by the NSF through the grants DMS-1162638 (Beck),
DMS-0946431 (Young), and DMS-1156499 (MSRI-UP REU), and by the NSA through grant H98230-11-1-0213.}

\maketitle

%------------------------------------------------%

% \section{Introduction}

Graph coloring problems are ubiquitous in many areas within and outside of mathematics.
We are interested in certain enumerative questions about coloring signed graphs.
A \emph{signed graph} $\Sigma = (\Gamma,\sigma)$ consists of a graph $\Gamma = (V,E)$ and a signature $\sigma \in \left\{ \pm
\right\}^E$.
The underlying graph $\Gamma$ may have multiple edges and, besides the usual links and loops, also \emph{half edges} (with only one
endpoint) and \emph{loose edges} (no endpoints); the last are irrelevant for coloring questions, and so we assume in this paper
that $\Sigma$ has no loose edges.
% The \emph{order} of $\Sigma$ is the number of nodes, written $n$.
An unsigned graph can be realized by a signed graph all of whose edges are labelled with $+$.
Signed graphs originated in the social sciences and have found applications also in biology, physics, computer science, and economics; see \cite{zaslavskydynamicsurvey} for a comprehensive bibliography.

% For a positive integer $n$, let $[n] := \left\{ 1, 2, \dots, n \right\}$ and $[\pm n] := \left\{ -n, -n+1, \dots, n \right\}$.
The \emph{chromatic polynomial} $c_\Sigma(2k+1)$ counts the \emph{proper $k$-colorings} $\x \in % [\pm k]^V$, 
\left\{ 0, \pm 1, \dots, \pm k \right\}^V $, namely, those colorings that satisfy for any edge $vw \in E$
\[
  x_v \ne \sigma_{vw} \, x_w 
\]
and $x_v \ne 0$ for any $v \in V$ incident with some half edge.
Zaslavsky \cite{zaslavskysignedcoloring} proved that $c_\Sigma(2k+1)$ is indeed a polynomial in $k$.
It comes with a companion, the \emph{zero-free chromatic polynomial} $c^*_\Sigma(2k)$, which counts all proper $k$-colorings $\x \in
% ([\pm k] \setminus 0)^V$.
\left\{ \pm 1, \dots, \pm k \right\}^V $.

\begin{figure}[htb]
\centering
\begin{tabular}{ l c r }
\begin{tikzpicture}[style=thin,]
\draw (18:1cm) -- (90:1cm) -- (162:1cm) -- (234:1cm) --
(306:1cm) -- cycle;
\draw (18:.5cm) -- (162:.5cm) -- (306:.5cm) -- (90:.5cm) --
(234:.5cm) -- cycle;
\foreach \x in {18,90,162,234,306}{
\draw (\x:.5cm) -- (\x:1cm);
\draw (\x:1cm) circle (2pt)[fill=black];
\draw (\x:.5cm) circle (2pt)[fill=black];
}
\node (A) at (270:1.2cm) {$P_1$}; 
\end{tikzpicture} & 
\begin{tikzpicture}[style=thin,]
%\draw (18:1cm) -- (90:1cm) -- (162:1cm) -- (234:1cm) --
(306:1cm) -- cycle;
\draw (18:1cm) -- (90:1cm) -- (162:1cm)--(234:1cm); 
\draw[dashed](234:1cm) --(306:1cm);
\draw(306:1cm)--(18:1cm);
\draw (18:.5cm) -- (162:.5cm) -- (306:.5cm) -- (90:.5cm) --
(234:.5cm) -- cycle;
\foreach \x in {18,90,162,234,306}{
\draw (\x:.5cm) -- (\x:1cm);
\draw (\x:1cm) circle (2pt)[fill=black];
\draw (\x:.5cm) circle (2pt)[fill=black];
} 
\node (A) at (270:1.2cm) {$P_2$};
\end{tikzpicture} & 
\begin{tikzpicture}[style=thin,]
\draw (18:1cm) -- (90:1cm);
\draw (90:1cm)--(162:1cm);
%\draw[dashed] (90:1cm)-- (162:1cm);
\draw[dashed] (162:1cm) -- (234:1cm);
\draw [dashed](306:1cm) -- (18:1cm);
\draw (234:1cm)--(306:1cm);
\draw (18:.5cm) -- (162:.5cm) -- (306:.5cm) -- (90:.5cm) --
(234:.5cm) -- cycle;
\foreach \x in {18,90,162,234,306}{
\draw (\x:.5cm) -- (\x:1cm);
\draw (\x:1cm) circle (2pt)[fill=black];
\draw (\x:.5cm) circle (2pt)[fill=black];
}
\node (A) at (270:1.2cm) {$P_3$};
\end{tikzpicture} \\
\begin{tikzpicture}[style=thin,]
\draw (18:1cm) -- (90:1cm);
\draw (90:1cm)--(162:1cm);
%\draw[dashed] (90:1cm)-- (162:1cm);
\draw (162:1cm) -- (234:1cm);
\draw (306:1cm) -- (18:1cm);
\draw [dashed](234:1cm)--(306:1cm);
\draw [dashed](18:.5cm) -- (162:.5cm) ;
\draw (162:.5cm) -- (306:.5cm) -- (90:.5cm) --
(234:.5cm)--(18:.5cm);
\foreach \x in {18,90,162,234,306}{
\draw (\x:.5cm) -- (\x:1cm);
\draw (\x:1cm) circle (2pt)[fill=black];
\draw (\x:.5cm) circle (2pt)[fill=black];
}
\node (A) at (270:1.2cm) {$P_4$};
\end{tikzpicture} &
\begin{tikzpicture}[style=thin,]
\draw (18:1cm) -- (90:1cm);
\draw (90:1cm)--(162:1cm);
%\draw[dashed] (90:1cm)-- (162:1cm);
\draw[dashed] (162:1cm) -- (234:1cm);
\draw [dashed](306:1cm) -- (18:1cm);
\draw (234:1cm)--(306:1cm);
\draw [dashed](18:.5cm) -- (162:.5cm) ;
\draw (162:.5cm) -- (306:.5cm) -- (90:.5cm) --
(234:.5cm)--(18:.5cm);
\foreach \x in {18,90,162,234,306}{
\draw (\x:.5cm) -- (\x:1cm);
\draw (\x:1cm) circle (2pt)[fill=black];
\draw (\x:.5cm) circle (2pt)[fill=black];
}
\node (A) at (270:1.2cm) {$P_5$};
\end{tikzpicture} &
\begin{tikzpicture}[style=thin,]
\draw (18:1cm) -- (90:1cm);
\draw (90:1cm)--(162:1cm);
%\draw[dashed] (90:1cm)-- (162:1cm);
\draw (162:1cm) -- (234:1cm);
\draw (306:1cm) -- (18:1cm);
\draw [dashed](234:1cm)--(306:1cm);
\draw [dashed](18:.5cm) -- (162:.5cm) ;
\draw (162:.5cm) -- (306:.5cm) -- (90:.5cm) --
(234:.5cm)--(18:.5cm);
\foreach \x in {18,162,234,306}{
\draw (\x:.5cm) -- (\x:1cm);
\draw [dashed] (90:.5cm)--(90:1cm);
\draw (90:1cm) circle (2pt)[fill=black];
\draw (90:.5cm) circle (2pt)[fill=black];
\draw (\x:1cm) circle (2pt)[fill=black];
\draw (\x:.5cm) circle (2pt)[fill=black];
}
\node (A) at (270:1.2cm) {$P_6$};
\end{tikzpicture} \\

\end{tabular}
\caption{The six switching-distinct signed Petersen graphs.}\label{sixsignedpetersen}
\end{figure}
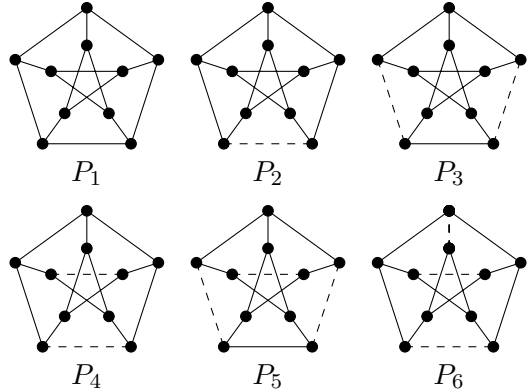

The \emph{Petersen graph} has served as a reference point to many
proposed results in graph theory. Considering \emph{signed}
Petersen graphs, Zaslavsky \cite{zaslavskypetersen} showed that, while there are $2^{15}$ ways to assign a signature to the
fifteen edges, only six of these are different up to switching isomorphism (a
notion that we will make precise below), % in Section \ref{background}), 
depicted in Figure~\ref{sixsignedpetersen}.
(In our figures we represent a positive edge with a solid line and a negative edge with a dashed line.)

In \cite{zaslavskypetersen} Zaslavsky proved that these six signed Petersen graphs have distinct chromatic polynomials; thus they can be distinguished by this signed-graph invariant. He did not compute the chromatic
polynomials but showed that they evaluate to distinct numbers at 3 \cite[Table~9.2]{zaslavskypetersen}.
He conjectured that the six different signed Petersen graphs also have distinct zero-free chromatic polynomials, and that both types of chromatic polynomials have distinct evaluations at \emph{any} positive integer
\cite[Conjecture~9.1]{zaslavskypetersen}.
Our first result confirms this conjecture.

\begin{theorem}
The chromatic polynomials of the signed Petersen graphs (denoted $P_1, \dots, P_6$ in Figure \ref{sixsignedpetersen}) are
\begin{align*}
  c_{P_1}(2k+1) &= 1024k^{10}-2560k^9+3840k^8-4480k^7+3712k^6 \\
                &\qquad -1792k^5+160k^4+480k^3-336k^2+72k \, , \\
  c_{P_2}(2k+1) &= 1024k^{10}-2560k^9+3840k^8-4480k^7+3968k^6 \\
                &\qquad -2560k^5+1184k^4-352k^3+48k^2 \, , \\
  c_{P_3}(2k+1) &= 1024k^{10}-2560k^9+3840k^8-4480k^7+4096k^6 \\
                &\qquad -2944k^5+1696k^4-760k^3+236k^2-40k \, , \\
  c_{P_4}(2k+1) &= 1024k^{10}-2560k^9+3840k^8-4480k^7+4224k^6 \\
                &\qquad -3200k^5+1984k^4-952k^3+308k^2-52k \, , \\
  c_{P_5}(2k+1) &= 1024k^{10}-2560k^9+3840k^8-4480k^7+4096k^6 \\
                &\qquad -3072k^5+1920k^4-960k^3+320k^2-48k \, , \\
  c_{P_6}(2k+1) &= 1024k^{10}-2560k^9+3840k^8-4480k^7+4480k^6 \\
                &\qquad -3712k^5+2560k^4-1320k^3+460k^2-90k \, .
\end{align*}
Their zero-free counterparts are
\begin{align*}
  c^*_{P_1}(2k) &= 1024k^{10} - 7680k^9 + 26880k^8 - 58240k^7 + 86592k^6  \\
&\qquad -91552k^5 + 68400k^4 - 34440k^3 + 10424k^2 - 1408k \, , \\
  c^*_{P_2}(2k) &= 1024k^{10}-7680k^9+26880k^8-58240k^7+86848k^6 \\
                &\qquad -93088k^5+72304k^4-39880k^3+14792k^2 - 3288k \, , \\
  c^*_{P_3}(2k) &= 1024k^{10} - 7680k^9 + 26880k^8 - 58240k^7 + 86976k^6 \\
                &\qquad - 93856k^5 + 74256k^4 - 42592k^3 + 16960k^2 - 4222k \, , \\
  c^*_{P_4}(2k) &= 1024k^{10} - 7680k^9 + 26880k^8 - 58240k^7 + 87104k^6 \\
                &\qquad - 94496k^5 + 75664k^4 - 44320k^3 + 18192k^2 - 4698k \, , \\
  c^*_{P_5}(2k) &= 1024k^{10} - 7680k^9 + 26880k^8 - 58240k^7 + 86976k^6 \\
                &\qquad - 93984k^5 + 74800k^4 - 43560k^3 + 17840k^2 - 4616k \, , \\
  c^*_{P_6}(2k) &= 1024k^{10} - 7680k^9 + 26880k^8 - 58240k^7 + 87360k^6 \\
                &\qquad - 95776k^5 + 78480k^4 - 47760k^3 + 20640k^2 - 5660k \, .
\end{align*}
Consequently (as a quick computation with a computer algebra system shows), none of the difference polynomials $c_{ P_m } (2k+1) - c_{ P_n } (2k+1)$ and $c^*_{ P_m } (2k) - c^*_{ P_n } (2k)$, with $m \ne n$, have a positive integer root.
\end{theorem}

To compute the above polynomials, we developed and executed a computer program (running in {\tt SAGE} \cite{sage}) that efficiently determines the number of proper $k$-colorings for any signed graph.
We include our code in the appendix and at {\tt math.sfsu.edu/beck/papers/signedpetersen.sage}.
% Section~\ref{background} contains the necessary background about signed graphs and a description of our implementation.

We also used our program to compute the chromatic polynomials of all signed complete graphs up to five vertices; up to switching isomorphism, there are two signed $K_3$'s, three signed $K_4$'s, and seven signed $K_5$'s. As with the signed Petersen graphs, the chromatic polynomials distinguish these signed complete graphs:

\begin{figure}[hbtp]
\begin{tabular}{ c c c c c }

\begin{tikzpicture}[style=thin,]
\draw(90:1cm) --(210:1cm) --(330:1cm) -- cycle;
\foreach \x in {90,210,330}{
\draw (\x:1cm) circle (2pt)[fill=black];
}
\node (A) at (270:1.2cm) {$K_3^{(1)}$};
\end{tikzpicture} 

\begin{tikzpicture}[style=thin,]
\draw  (330:1cm)--(90:1cm) --(210:1cm);
\draw [dashed] (210:1cm)-- (330:1cm);
\foreach \x in {90,210,330}{
\draw (\x:1cm) circle (2pt)[fill=black];
}
\node (A) at (270:1.2cm) {$K_3^{(2)}$};
\end{tikzpicture} 
\begin{tikzpicture}[style=thin,]
\draw (45:1cm)--(135:1cm)--(225:1cm)--(315:1cm)--cycle;
\draw (45:1cm)--(225:1cm);
\draw (135:1cm)--(315:1cm);
\foreach \x in {45,135,225,315}{
\draw (\x:1cm) circle (2pt)[fill=black];
}
\node (A) at (270:1.2cm) {$K_4^{(1)}$};
\end{tikzpicture} 
\begin{tikzpicture}[style=thin,]
\draw (315:1cm)--(45:1cm)--(135:1cm)--(225:1cm);
\draw[dashed](225:1cm)--(315:1cm);
\draw (45:1cm)--(225:1cm);
\draw (135:1cm)--(315:1cm);
\foreach \x in {45,135,225,315}{
\draw (\x:1cm) circle (2pt)[fill=black];
}
\node (A) at (270:1.2cm) {$K_4^{(2)}$};
\end{tikzpicture} 
\begin{tikzpicture}[style=thin,]
\draw (315:1cm)--(45:1cm);
\draw (135:1cm)--(225:1cm);
\draw [dashed] (45:1cm) --(135:1cm);
\draw[dashed](225:1cm)--(315:1cm);
\draw (45:1cm)--(225:1cm);
\draw (135:1cm)--(315:1cm);
\foreach \x in {45,135,225,315}{
\draw (\x:1cm) circle (2pt)[fill=black];
}
\node (A) at (270:1.2cm) {$K_4^{(3)}$};
\end{tikzpicture} 

\end{tabular} \\

\begin{tabular}{ c c c c}
\begin{tikzpicture}[style=thin,]
\draw (90:1) -- (162:1) -- (234:1) -- (306:1) --(18:1)--cycle;
\draw (234:1)--(90:1) -- (306:1) --(162:1) -- (18:1)--cycle;
\foreach \x in {90,162,234,306,18}{
\draw (\x:1cm) circle (2pt)[fill=black];
}
\node (A) at (270:1.2cm) {$K_5^{(1)}$};
\end{tikzpicture}
\begin{tikzpicture}[style=thin,]
\draw (90:1) -- (162:1) -- (234:1);
\draw[dashed](234:1) --(306:1);
\draw(306:1)--(18:1)--(90:1);
\draw (234:1)--(90:1) -- (306:1) --(162:1) -- (18:1)--cycle;
\foreach \x in {90,162,234,306,18}{
\draw (\x:1cm) circle (2pt)[fill=black];
}
\node (A) at (270:1.2cm) {$K_5^{(2)}$};
\end{tikzpicture}
\begin{tikzpicture}[style=thin,]
\draw (90:1) -- (162:1) -- (234:1);
\draw[dashed](234:1) --(306:1);
\draw(306:1)--(18:1)--(90:1);
\draw (234:1)--(90:1) -- (306:1)--(162:1);
\draw[dashed](162:1) -- (18:1);
\draw (18:1)--(234:1);
\foreach \x in {90,162,234,306,18}{
\draw (\x:1cm) circle (2pt)[fill=black];
}
\node (A) at (270:1.2cm) {$K_5^{(3)}$};
\end{tikzpicture}
\begin{tikzpicture}[style=thin,]
\draw [dashed] (18:1)--(90:1)--(162:1);
\draw (162:1) -- (234:1) -- (306:1) --(18:1);
\draw (234:1)--(90:1) -- (306:1) --(162:1) -- (18:1)--cycle;
\foreach \x in {90,162,234,306,18}{
\draw (\x:1cm) circle (2pt)[fill=black];
}
\node (A) at (270:1.2cm) {$K_5^{(4)}$};
\end{tikzpicture}

\end{tabular} \\
\begin{tabular} {c c c}
\begin{tikzpicture}[style=thin,]
\draw [dashed] (18:1)--(90:1)--(162:1);
\draw[dashed] (234:1)--(306:1);
\draw (162:1) -- (234:1);
\draw (306:1) --(18:1);
\draw (234:1)--(90:1) -- (306:1) --(162:1) -- (18:1)--cycle;
\foreach \x in {90,162,234,306,18}{
\draw (\x:1cm) circle (2pt)[fill=black];
}
\node (A) at (270:1.2cm) {$K_5^{(5)}$};
\end{tikzpicture}
\begin{tikzpicture}[style=thin,]
\draw (18:1)--(90:1)--(162:1);
\draw[dashed]  (162:1) -- (234:1)--(306:1) --(18:1);
\draw (234:1)--(90:1) -- (306:1) --(162:1) -- (18:1)--cycle;
\foreach \x in {90,162,234,306,18}{
\draw (\x:1cm) circle (2pt)[fill=black];
}
\node (A) at (270:1.2cm) {$K_5^{(6)}$};
\end{tikzpicture}
\begin{tikzpicture}[style=thin,]
\draw [dashed] (18:1)--(90:1)--(162:1);
\draw[dashed] (234:1)--(306:1);
\draw (162:1) -- (234:1);
\draw (306:1) --(18:1);
\draw (234:1)--(90:1) -- (306:1)--(162:1);
\draw[dashed](162:1) -- (18:1);
\draw (18:1)--(234:1);
\foreach \x in {90,162,234,306,18}{
\draw (\x:1cm) circle (2pt)[fill=black];
}
\node (A) at (270:1.2cm) {$K_5^{(7)}$};
\end{tikzpicture}
\end{tabular}
\caption{The switching classes of signed complete graphs.}\label{signedcomplete}
\end{figure}
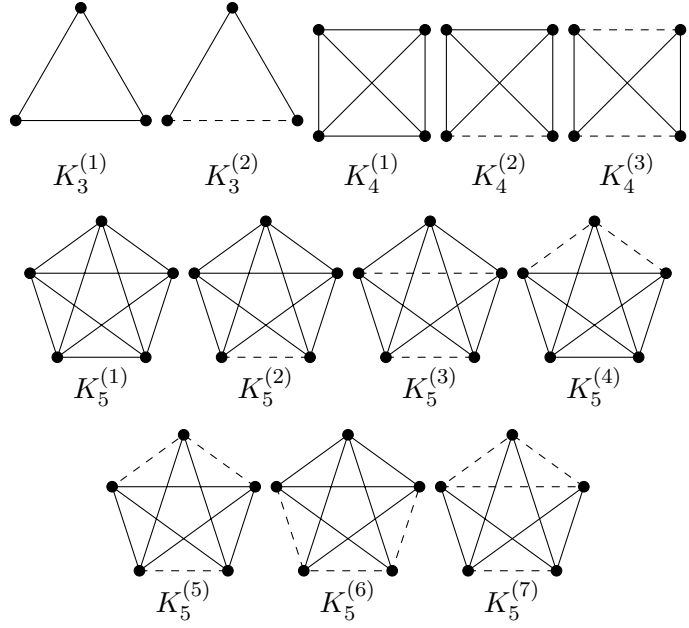

\begin{theorem}
The chromatic polynomials of the signed complete graphs (denoted $K_3^{ (1) }, K_3^{ (2) }, \dots, K_5^{ (7) }$ in Figure
\ref{signedcomplete}) are
\begin{align*}
  c_{K_3^{ (1) }}(2k+1) &= 8k^3 - 2k \, , \\
  c_{K_3^{ (2) }}(2k+1) &= 8k^3 \, , \\
  c_{K_4^{ (1) }}(2k+1) &= 16k^4-16k^3-4k^2+4k \, , \\
  c_{K_4^{ (2) }}(2k+1) &= 16k^4-16k^3+4k^2 \, , \\
  c_{K_4^{ (3) }}(2k+1) &= 16k^4-16k^3+12k^2-2k \, , \\
 c_{K_5^{ (1) }}(2k+1) &= 32k^5 - 80k^4 + 40k^3 + 20k^2 - 12k \, , \\
c_{K_5^{ (2) }}(2k+1) &= 32k^5 - 80k^4 + 64k^3 - 16k^2 \, , \\
c_{K_5^{ (3) }}(2k+1) &= 32k^5 - 80k^4 + 88k^3 - 48k^2 + 10k \, , \\ 
c_{K_5^{ (4) }}(2k+1) &= 32k^5 - 80k^4 + 72k^3 - 28k^2 + 4k \, . \\
c_{K_5^{ (5) }}(2k+1) &= 32k^5 - 80k^4 + 96k^3 - 56k^2 + 12k \, , \\
c_{K_5^{ (6) }}(2k+1) &= 32k^5 - 80k^4 + 80k^3 - 40k^2 + 8k \, , \\
c_{K_5^{ (7) }}(2k+1) &= 32k^5 - 80k^4 + 120k^3 - 80k^2 + 20k \, .
\end{align*}
\end{theorem}

The corresponding zero-free chromatic polynomials are
\begin{align*}
 c^*_{K_3^{ (1) }}(2k) &= 8k^3 - 12k^2+4k \, , \\
 c^*_{K_3^{ (2) }}(2k) &= 8k^3-12k^2+6k \, , \\
 c^*_{K_4^{ (1) }}(2k) &= 16k^4-48k^3+44k^2-12k \, , \\
  c^*_{K_4^{ (2) }}(2k) &= 16k^4-48k^3+52k^2-24k \, , \\
  c^*_{K_4^{ (3) }}(2k) &= 16k^4-48k^3+60k^2-34k \, , \\
 c^*_{K_5^{ (1) }}(2k) &= 32k^5 - 160k^4 + 280k^3 - 200k^2 +48k \, , \\
c^*_{K_5^{ (2) }}(2k) &= 32k^5 - 160k^4 + 304k^3 - 272k^2 +114k \, , \\
c^*_{K_5^{ (3) }}(2k) &= 32k^5 - 160k^4 + 328k^3 - 340k^2 +174k \, , \\
c^*_{K_5^{ (4) }}(2k) &= 32k^5 - 160k^4 + 312k^3 - 296k^2 +136k \, , \\
c^*_{K_5^{ (5) }}(2k) &= 32k^5 - 160k^4 +336k^3 - 360k^2 +190k \, , \\
c^*_{K_5^{ (6) }}(2k) &= 32k^5 - 160k^4 + 320k^3 - 320k^2 +158k \, , \\
c^*_{K_5^{ (7) }}(2k) &= 32k^5 - 160k^4 + 360k^3 - 420k^2 +240k \, .
\end{align*}

%------------------------------------------------%

% \section{Background and Implementation}\label{background}

We now review a few constructs on a signed graph $\Sigma = (V, E, \sigma)$ and describe our implementation.
The \emph{restriction} of $\Sigma$ to an edge set $F \subseteq E$ is the signed graph $\left( V, F, \sigma|_F \right)$.
For $e \in E$, we denote by $\Sigma - e$ (the \emph{deletion} of $e$) the restriction of $\Sigma$ to $E - \{ e \}$.
For $v \in V$, denote by $\Sigma - v$ the restriction of $\Sigma$ to $E-F$ where $F$ is the set of all edges incident to~$v$.
A component of the signed graph $\Sigma = (\Gamma, \sigma)$ is \emph{balanced} if it contains no half edges and each cycle has positive sign product.

\emph{Switching} $\Sigma$ by $s \in \left\{ \pm \right\}^V$ results in the new signed graph $(V, E, \sigma^s)$ where
$
  \sigma^s_{ vw } = s_v \, \sigma_{ vw } \, s_w
$.
Switching does not alter balance, and any balanced signed graph can be obtained from switching an all-positive graph
\cite{zaslavskysignedgraphs}.
We also note that there is a natural bijection of proper colorings of $\Sigma$ and a switched version of it, and this bijection preserves
the number of proper $k$-colorings.
Thus the chromatic polynomials of $\Sigma$ are invariant under switching.

The \emph{contraction} of $\Sigma$ by $F \subseteq E$, denoted by $\Sigma/F$, is defined as follows \cite{zaslavskysignedgraphs}:
switch $\Sigma$ so that every balanced component of $F$ is all positive, coalesce all nodes of each balanced component, and discard
the remaining nodes and all edges in $F$; note that this may produce half edges.
If $F = \left\{ e \right\}$ for a link $e$, $\Sigma/e$ is obtained by switching $\Sigma$ so that $\sigma(e) = +$ and
then contracting $e$ as in the case of unsigned graphs, that is, disregard $e$ and identify its two endpoints.
If $e$ is a negative loop at $v$, then $\Sigma/e$ has vertex set $V - \{ v \}$ and edge set resulting from $E$
by deleting $e$ and converting all edges incident with $v$ to half edges.
The chromatic polynomial satisfies the deletion--contraction formula \cite{zaslavskysignedcoloring}
\begin{equation}\label{deletioncontraction}
  c_\Sigma(2k+1) = c_{\Sigma - e} (2k+1) - c_{\Sigma/e} (2k+1) \, .
\end{equation}
The zero-free chromatic polynomial $c_\Sigma^*(2k)$ satisfies the same identity provided that $e$ is not a half edge or negative loop.
We will use \eqref{deletioncontraction} repeatedly in our computations.

We encode a signed graph $\Sigma$ by its \emph{incidence matrix} as follows: 
first \emph{bidirect} $\Sigma$, i.e., give each edge an independent orientation at each endpoint (which we think of as an arrow pointing towards or
away from the endpoint), such that a positive edge has one arrow pointing towards one and away from the other endpoint, and a negative edge has
both arrows pointing either towards or away from the endpoints.
The incidence matrix has rows indexed by vertices, columns indexed by edges, and entries equal to $\pm 1$ according to whether the edge points towards or
away from the vertex (and 0 otherwise).
Since half edges and negative loops have the same effect on the chromatic polynomial of $\Sigma$, we may assume that $\Sigma$ has no half edge.

\begin{figure}
\begin{tabular} { l  r }
\tikzstyle{every node}=[circle, draw, fill=black, inner sep=0pt, minimum width=4pt]
\begin{tikzpicture}[style=thin,]
\draw (0,0)--(0,2);
\draw (2,0)--(2,2);
\draw [dashed] (0,0)--(2,0);
\draw [dashed] (0,2)--(2,2);
\draw (0,0)--(2,2);
\draw (2,0)--(0,2);
\draw (0,0) circle (2pt)[fill=black];
\draw (0,2) circle (2pt)[fill=black];
\draw (2,0)  circle (2pt)[fill=black];
\draw (2,2) circle (2pt)[fill=black];
\draw (2.25,-.25) node[draw=none,fill=none] {$c$};
\draw (2.25,2.25) node[draw=none,fill=none] {$b$};
\draw (-.25,-.25) node[draw=none,fill=none] {$d$};
\draw (-.25,2.25) node[draw=none,fill=none] {$a$};
\end{tikzpicture}

\tikzstyle{every node}=[circle, draw, fill=black, inner sep=0pt, minimum width=4pt]
\begin{tikzpicture}[style=thin,]
\draw (0,0)--(0,2);
\draw (2,0)--(2,2);
\draw (0,0)--(2,0);
\draw (0,2)--(2,2);
\draw (0,0)--(2,2);
\draw (2,0)--(0,2);
\draw [->] (0,0)--(0,.75);
\draw [->] (0,1)--(0,1.25);
\draw [->] (2,0)--(2,.75);
\draw [->] (2,1)--(2,1.25);
\draw [->] (.5,.5)--(.75,.75);
\draw [->] (1,1)--(1.25,1.25);
\draw [->] (.5,1.5)--(.75,1.25);
\draw [->] (1,1)--(1.25,.75);
\draw [->] (0,0)--(.75,0);
\draw [<-] (1,0)--(1.25,0);
\draw [->] (0,2)--(.75,2);
\draw [<-] (1,2)--(1.25,2);
\draw (0,0) circle (2pt)[fill=black];
\draw (0,2) circle (2pt)[fill=black];
\draw (2,0)  circle (2pt)[fill=black];
\draw (2,2) circle (2pt)[fill=black];
\draw (2.25,-.25) node[draw=none,fill=none] {$c$};
\draw (2.25,2.25) node[draw=none,fill=none] {$b$};
\draw (-.25,-.25) node[draw=none,fill=none] {$d$};
\draw (-.25,2.25) node[draw=none,fill=none] {$a$};
\end{tikzpicture}
\end{tabular}

$\begin{matrix} & ab & ac & ad & bc & bd & cd  \\
a & -1 &-1 &  1 & 0 & 0 & 0\\
b & -1 &0 &  0 & 1 & 1 & 0\\
c & 0 &1 &  0 & -1 & 0 & -1\\
d & 0 &0 &  -1 & 0 & -1 & -1 \\ \end{matrix}$
\caption{ $K_4^{(3)}$ with one of its bidirections and corresponding incidence matrix.} \label{signedK4}
\end{figure}
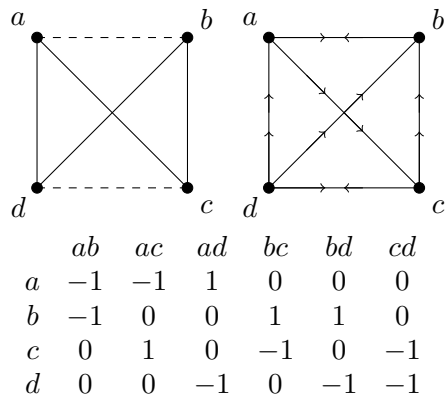 

Deletion--contraction can be easily managed by incidence matrices:
deletion of an edge simply means deletion of the corresponding column;
contraction of a positive edge $vw$ means replacing the rows corresponding to $v$ and $w$ by their sum and then deleting the column corresponding to the edge $vw$ (it is
sufficient to only consider contraction of positive edges, since we can always switch one of its endpoints if necessary, which means negating the
corresponding row). 
Note that this process works for both links and half edges.
Note also that we will constantly look for multiple edges (with the same sign) and replace them with a single edge.

Thus we can keep track of incidence matrices as we recursively apply deletion--contraction, leading to empty signed graphs or signed graphs that only have half edges; both
have easy chromatic polynomials.

%------------------------------------------------%

\section*{Appendix: Code}

Below is the {\tt SAGE} code (which can be loaded into any {\tt SAGE} terminal) used to compute chromatic polynomials of signed graphs. The procededure \texttt{chrom} is the main method which takes an incidence matrix and outputs the chromatic polynomial as an expression.

\vspace{.3in}
\small
\begin{verbatim}

#All positive edges in columns have a 1 and -1 entry. Standardize makes any 
#positive edge column have the 1 come before the -1 when reading the column 
#top to bottom. All negative edges in columns have two 1 entries or two -1 entries. 
#Standardize makes any negative edge column have two 1s. Negative loops have
# one 1 or -1 entry in its edge column. Standardize makes all negative loop 
#columns have one 1. Standardizing our edges into one convention makes finding a  
#positiveedge to delete-contract more efficient. 
def standardize(G):
    G_rows=G.nrows()
    G_cols=G.ncols()
    for j in range(0,G_cols):
        for i in range(0,G_rows):
            if G[i,j] == 1:
                break
            if G[i,j] == -1:
                G.add_multiple_of_column(j,j,-2)
                break
    return G

#If there is a muti-edge in the peterson graph, meaning a column is repeated, then 
#delete the column and return a new incidence matrix.
def check(Gin):
    G=standardize(Gin) #
    G_cols=G.ncols() 
    G_rows=G.nrows()
    l=G.columns()
    for i in range(0,G_cols-1):
        for j in range(i+1,G_cols):
            if G.column(G_cols-1-i)== G.column(G_cols-1-j):
                l.pop(G_cols-1-i)
                break
    B=matrix(l)
    C=B.transpose()
    return C

#Returns an ordered pair of the incidence matrix of the graph with an edge
# deleted and the incidence matrix of the graph with an edge contracted. 
def DC(G):
    rows=G.nrows()
    cols=G.ncols()
    #If a single edge, return an empty graph with the same number of vertices for 
    #the graph with the edge deleted. Return an empty graph with one 
    #less vertex for the graph with the edge contracted.
    if cols == 1:           
        C=matrix(QQ,rows,1,range(0))
        H=matrix(QQ,rows-1,1,range(0))
        return (C,H)
    #Else, find the first positive edge in the matrix when reading the matrix
    # from left to right. Delete and contract this edge returning new incidence matrices.
    else:
        j=0
        #Increment through the columsn looking for a positive edge
        while j<cols:
            sum=0
            for i in range(0,rows):
                sum = G[i,j] + sum
        
            #If the non-sero entires in the column add to zero, this is a positive edge. 
             #Delete it.
            if sum == 0:
                    delete_col=j
                    break
            j=j+1

        #If we found a positive edge then j will not be the number of columns in the incidence 
        #matrix and we have deleted a column from the matrix. Create a new incidence matrix 
         #without this column sotred in D, D represents a new graph with an edge deleted.
        if j != cols:
            l=G.columns()
            l.pop(delete_col)
            B=matrix(l)
            C=B.transpose()
            D=B.transpose()
            #Use this incidence matrix to contract an edge. 
             #This is done by adding together the two 
             #rows which contained the deleted  column's non-zero entries 
             #creating a new incidence matrix H2.
            r1found = false
            for k in range(0,rows):
                if G[k,delete_col]!=0:
                    if r1found == false:
                        r_1=k
                        r1found = true
                    else:
                        r_2=k
	    C.add_multiple_of_row(r_2,r_1,1)
            for p in range(0,cols-1):
                if abs(C[r_2,p])==2:
                    C[r_2,p]=1
            F=C.rows()
            F.pop(r_1)
            H=matrix(F)
            H1=check(H)
            D1=standardize(D)
            H2=standardize(H1)
            return (D1,H2)
           
        #If there was no positive edge in the matrix, then there exists a negative loop. 
        #Delete Contract the negative loop.        
        if j==cols:
            return c_neg_loop(G)

#Returns an ordered pair of the incidence matrix of the graph with a negative loop deleted 
#and the incidence matrix of the graph with a negative loop contracted. 
def c_neg_loop(G):
    G_rows=G.nrows()
    G_cols=G.ncols()
    E=G
    for j in range(0,G_cols):
        sum = 0
        for i in range(0,G_rows):
            sum=G[i,j] + sum
            if abs(G[i,j]) == 1:
                r1=i
            
        if abs(sum) == 1:
            l=G.columns()
            l.pop(j)
            B=matrix(l)
            D=B.transpose()
            m=D.rows()
            m.pop(r1)
            E=matrix(m)
            D1=standardize(D)
            E1=standardize(E)
            return (D1,E1)
            break       
    if abs(sum)!= 1:
        return "All negative edges and no loops!"

#To ensure there is a positive edge or negative loop to delete-contract, check 
#to make sure one exists in the incidence matrix. If not, we can switch a vertex 
#by negating a row, creating a positive edge or negative loop to delete and contract.
def switch(M):
    posEdge=False
    for item in M.columns():
        vNum=0
        sum=0
        nzExists=false
        for i in range (0,M.nrows()):
            if item[i] !=0:
                nzExists=True
                rowNum=i
            sum=sum+item[i]
        if nzExists==True and sum==0:
            posEdge=True
        elif abs(sum)==1:
	    posEdge=True
        elif nzExists==True and abs(sum)==2:
            vNum=rowNum
    if posEdge==False and M.ncols()!=0:
        M.add_multiple_of_row(vNum,vNum,-2)
    return M

#Returns the chromatic polynomial of a graph represented as an incidence matrix.    
def chrom(P):
    #Make sure the graph has a positive edge to delete and contract and delete any 
     #multiple edges using switch and check.
    P=switch(P)
    P=check(P)
    #If the graph is a negative loop, return the chromatic polynomial.
    if len(P.columns())==1:
        if len(P.rows())==1:
		if P[0][0]==abs(1):
			return 2*x
    #If the graph is empty, return the chromatic polynomial.
    empty=true
    for item in P.columns()[0]:
	if item!=0:
		empty=false
    if empty==true:
	return ((2*x)+1)^(P.nrows())
    #Else, delete and contract the graph recursively calling chrom on both the
    # graph with the deleted edge and the graph with the contracted edge. 
     #Return the difference of the parts.               
    (z,y)=DC(P)
    return expand(chrom(z))-expand(chrom(y))
    
#Returns the chromatic polynomial of a graph represented as an incidence matrix.    
def zeroFreeChrom(P):
    #Make sure the graph has a positive edge to delete and contract and delete any 
     #multiple edges using switch and check.
    P=switch(P)
    P=check(P)
    #If the graph is a negative loop, return the chromatic polynomial.
    if len(P.columns())==1:
        if len(P.rows())==1:
		if P[0][0]==abs(1):
			return 2*x
    #If the graph is empty, return the chromatic polynomial. 
     #There are 2x choices for every vertex in 
    #the graph because we are not coloring with zero
    empty=true
    for item in P.columns()[0]:
	if item!=0:
		empty=false
    if empty==true:
	return ((2*x))^(P.nrows())
    #Else, delete and contract the graph recursively calling zeroFreeChrom on 
     #both the graph with the deleted edge and the graph with the contracted edge. 
      #Return the difference of the parts.               
    (z,y)=DC(P)
    return expand(zeroFreeChrom(z))-expand(zeroFreeChrom(y))
\end{verbatim}
\normalsize

\vspace{.3in}

%------------------------------------------------%

\bibliographystyle{amsplain}
%\bibliography{signedpetersen} 

\begin{thebibliography}{1}

\bibitem{sage}
William A.~Stein et~al, \emph{Sage mathematics software, version 5.1} (2012),
  {\tt http://www.sagemath.org}.

\bibitem{zaslavskysignedcoloring}
Thomas Zaslavsky, \emph{Signed graph coloring}, Discrete Math. \textbf{39}
  (1982), no.~2, 215--228.

\bibitem{zaslavskysignedgraphs}
\bysame, \emph{Signed graphs}, Discrete Appl. Math. \textbf{4} (1982),
  47--74. Erratum, Discrete Appl. Math. \textbf{5} (1983), 248.

\bibitem{zaslavskydynamicsurvey}
\bysame, \emph{A mathematical bibliography of signed and gain graphs and allied
  areas}, Electron. J. Combin. \textbf{5} (1998), Dynamic Surveys 8, 124 pp.
  (electronic), Electronically available at {\tt
  http://www.math.binghamton.edu/zaslav/Bsg/index.html}.

\bibitem{zaslavskypetersen}
\bysame, \emph{Six signed {P}etersen graphs, and their automorphisms}, Discrete
  Math. \textbf{312} (2012), no.~9, 1558--1583.

\end{thebibliography}

\providecommand{\bysame}{\leavevmode\hbox to3em{\hrulefill}\thinspace}
\providecommand{\MR}{\relax\ifhmode\unskip\space\fi MR }
\providecommand{\MRhref}[2]{%
  \href{http://www.ams.org/mathscinet-getitem?mr=#1}{#2}
}
\providecommand{\href}[2]{#2}

\end{document}